\documentclass[12pt,a4paper]{article}
\usepackage{amsfonts}
\usepackage{amsmath}
\usepackage{amssymb}
\usepackage[left=3cm,right=2cm,top=2cm,bottom=2cm]{geometry}

\begin{document}
\large
\centerline{\textbf{ On matrices with simple spectra arising from tensor products}}
\vskip0.5cm
\centerline{R.N.Gumerov, S.I.Vidunov }
\vskip0.5cm
\centerline{\textit{Kazan Federal University}}
\vskip0.5cm
\centerline{\textit{e-mail:renat.gumerov@kpfu.ru}}
\vskip1cm

\centerline{\textbf{Introduction.}}

\vskip0.5cm

Throughout all matrices are square with the real or complex entries.
The~purpose of this note is an application of topological methods to
a problem on simultaneous approximation of several matrices by a finite family of diagonalizable matrices satisfying
an additional spectral condition.  That is the simplicity
for all eigenvalues of a matrix product.

Problems of approximating  tuples of matrices and linear operators by finite collections of those elements with some extra nice properties have a considerable history.
They are intimately related to the very interesting questions in linear algebra, functional analysis, algebraic geometry and topology (see, e.g., [O'MCV], [DS] and the references therein).
We mention the result established in 1955
by Motzkin and Taussky ([MT], Theorem~5) who proved that any pair of complex commuting
matrices is approximately simultaneously diagonalizable.
Moreover, those problems  continue to find new
applications, e.g.,
to the study of phylogenetic invariants in biomathematics
[AR],[O'MCV, Ch.~6].

Originally we were motivated by our interest ([GG],[G1],[G2]) in approximating
elements in  the tensor products of Banach spaces
and  applications of parameterized families of polynomials
to topological groups, particularly
to perturbations in the matrix groups. Matrix tuples and  products of their elements
arise naturally in several contexts in the study of tensors, or multi-dimensional arrays,
especially in connection with bounds on tensor ranks (see [AS],[AL],[TB],[S],[SMS],[T]). A part of the
motivation for our work comes from the article [T] where some nontrivial tensor ranks estimations
for inverse matrices are presented.

\vskip1cm

\centerline{\textbf{Preliminaries.}}

\vskip0.5cm

Let $\mathbb{F}$ denote either the field of complex numbers $\mathbb{C}$ or the field of real numbers $\mathbb{R}$. The linear space of all square $n\times n$ matrices $M_n(\mathbb{F})$ over the field $\mathbb{F}$ is endowed with an arbitrary norm $\|\cdot\|$ which generates the same
metric topology .
The space $M_n(\mathbb{F})$ with
the ordinary matrix multiplication is a Banach algebra. As usual, the set
of all
invertible matrices in $M_n(\mathbb{F})$ is denoted by  $GL_n(\mathbb{F})$.
It is called \textit{the general linear group of degree $n$.}
\emph{An eigenvalue} of a matrix is said to be \emph{simple} if its algebraic
multiplicity equals one.
\emph{The spectrum} of a matrix is said to be \emph{simple} provided that all
eigenvalues of a given matrix are simple. In other words, if all eigenvalues are
pairwise distinct. Sometimes the matrices with simple spectra
are called \emph{generic}. Certainly, such matrices are diagonalizable.

In the following theorem we list necessary well known facts about topological
properties of matrix sets.
\vskip0.5cm
\textbf{Theorem.} \emph{The following properties are fulfilled in the space $M_n(\mathbb{F})$: }
\begin{itemize}
\item \emph{The set  $GL_n(\mathbb{F})$ is  dense and open.}

\item \emph{The general linear group is a locally compact group.}


\item \emph{The set of all invertible matrices with simple spectra
is dense and open.}

\end{itemize}

\vskip0.5cm

Recall that a mapping $f:X \rightarrow Y$ between two topological spaces is said to be
\textit{open}, if for any open set $O$ in $X$ the image $f(O)$ is open in $Y$. If
$f:GL_n(\mathbb{F})\rightarrow GL_n(\mathbb{F})$ is a surjective continuous
homomorphism then it is open ([HR], Theorem~5.29).

\vskip0.5cm
\vskip0.5cm
\centerline{\textbf{Results.}}

\vskip0.5cm
In this section we prove several propositions concerning perturbations of matrix tuples by finite
collections of matrices with simple spectra.
To do this we make use of the facts from Theorem.
\vskip0.5cm

\textbf{Proposition~1.} \emph{Let $f$ and $g$  be self-mappings of
the general linear group of degree $n$ over the field $\mathbb{F}$. Suppose that at least one of those mappings is open.
Let $A$ and $B$ be  matrices in $M_n(\mathbb{F})$. Then for any $\varepsilon > 0$ there exists a pair of invertible matrices
$A_{\varepsilon}$ and $B_{\varepsilon}$ with simple spectra such that
$$
\|A-A_{\varepsilon}\|<\varepsilon, \, \|B-B_{\varepsilon}\|<\varepsilon
$$
and the product matrix  $f(A_{\varepsilon})g(B_{\varepsilon})$ has a simple spectrum.}
\vskip0.5cm

\textbf{Proof.}  For definiteness, we assume that the mapping
 $f:GL_n(\mathbb{F})\rightarrow GL_n(\mathbb{F})$ is open.
Fix an arbitrary real number $\varepsilon > 0.$

Using Theorem,
we choose two invertible matrices $\tilde{A}$ and $B_{\varepsilon}$ with simple spectra
satisfying the inequalities
$$
\|A-\tilde{A}\|<\frac{\varepsilon}{2}, \qquad \|B-B_{\varepsilon}\|<\varepsilon.
$$

Further we consider the right multiplication by the element $g(B_{\varepsilon})$:
$$\
R:GL_n(\mathbb{F})\longrightarrow GL_n(\mathbb{F}):X\longmapsto  Xg(B_{\varepsilon})
$$
which is an open mapping.

By Theorem, in the space $M_n(\mathbb{F})$ we can take  the open ball $B(\tilde{A};r)$ of radius $r>0$ centered at the point $\tilde{A}$ , i.e., the open set
$$
B(\tilde{A};r)=\{X\in M_n(\mathbb{F}): \,\|\tilde{A}-X\|<r\}
$$
such that the following two conditions are satisfied:

\begin{itemize}
\item $r<\frac{\varepsilon}{2}$;
\item $B(\tilde{A};r)$ consists of invertible matrices with simple spectra.
\end{itemize}

Since the composition of open mappings is open the image
$R(f(B(\tilde{A};r)))$ is a non-empty open set.
It follows from  Theorem that there exists an invertible matrix
$C\in R(f(B(\tilde{A};r)))$ whose spectrum is simple.
Therefore the equality $C=f(A_{\varepsilon})g(B_{\varepsilon})$ holds
for some matrix $A_{\varepsilon}\in B(\tilde{A};r).$

Obviously, the matrices $A_{\varepsilon}$ and $B_{\varepsilon}$ are the required
perturbations of the given matrices $A$ and $B$. The proposition is proved.

\vskip0.5cm

Further, let  $f$ be the identity mapping and $g$ be the inverse mapping
on the group $GL_n(\mathbb{F})$.  We have the following consequence of Proposition~1.

\vskip0.5cm
\textbf{Proposition~2.} \emph{Let $A$ and $B$ be elements of $M_n(\mathbb{F})$. Then for any $\varepsilon > 0$ there exists a pair of invertible matrices
$A_{\varepsilon}$ and $B_{\varepsilon}$ with simple spectra such that
$$
\|A-A_{\varepsilon}\|<\varepsilon, \, \|B-B_{\varepsilon}\|<\varepsilon
$$
and the product matrix $A_{\varepsilon}B_{\varepsilon}^{-1}$ has a simple spectrum .}
\vskip0.5cm
There is an interesting  application of Proposition~2 to estimating the tensor rank
of the inverse matrix in [T] when the matrices
$A$ and $B$ are the factors of the Kronecker (tensor) products. More precisely, let
$X\in GL_{pq}(\mathbb{F})$ be a matrix whose tensor rank $tRank_{p,q}(X)=2$.
That is,  a decomposition of $X$ into the sum
of the Kronecker products of two matrices from  $M_p(\mathbb{F})$ and $M_q(\mathbb{F})$
with a minimal possible number of summands consists of two terms, i.e.,
$
X=A\otimes C + B\otimes D \quad \text{for some} \quad A,B\in M_p(\mathbb{F})\quad \text{and} \quad C,D\in M_q(\mathbb{F}).
$
Proposition~2 guarantees
the possibility of approximating the matrix $X$ by tensor-product binomials with
 two factors satisfying the properties indicated in that assertion. This fact is involved
 in the proof of the following upper bound (see [T,p.3172]):
 $$
 tRank_{p,q}(X^{-1})\leqslant \min\{p,q\}.
$$

The matrix products of the forms $AB^{-1}$ and $A^{-1}B$ with different spectral
conditions arise in studying the tensor ranks of $3$-tensors
(see, e.g., [AS],[AL],[TB],[S],[SMS]). In this case the slices of three-dimensional
arrays serve as the matrices $A$ and $B$.

It is worth noting that the spectral properties of such matrix products
play a significant role in a number of problems in the theory of matrix
pencils(see, e.g.,[I]).




In conclusion, we formulate the following statement which is a natural generalization of Proposition~1.
\vskip0.5cm

\textbf{Proposition~3.} \emph{Let $f_1,f_2,\ldots, f_k$   be a $k$-tuple of self-mappings
for
the general linear group of degree $n$ over the field $\mathbb{F}$. Suppose that at
least one of those mappings is open.
Let $A_1,A_2,\ldots, A_k$ be  a $k$-tuple of matrices in $M_n(\mathbb{F})$. Then for any $\varepsilon > 0$ there exists a $k$-tuple $A_{1\varepsilon},A_{2\varepsilon},\ldots, A_{k\varepsilon}$ consisting of  invertible matrices with simple spectra
 such that
$$
\|A-A_{1\varepsilon}\|<\varepsilon, \, \|A-A_{2\varepsilon}\|<\varepsilon, \ldots,\|A-A_{k\varepsilon}\|<\varepsilon,
$$
and the product matrix
$f_1(A_{1\varepsilon})f_2(A_{2\varepsilon})\ldots f_k(A_{k\varepsilon})$
has a simple spectrum.}
\vskip0.5cm
\centerline{ACKNOWLEDGMENTS}

It is a pleasure to thank our colleagues and friends for helpful conversations about
the results in this note. We are grateful to the participants of the seminars at KFU
and  KSPEU for their interest in the subject.

\end{document}